\documentclass[twoside, 11pt]{article}

\usepackage{amsmath}
\usepackage{paralist}
\usepackage{amssymb,mathrsfs,psfrag,graphicx}
\usepackage{hyperref}
\usepackage{subcaption}
\def\x{\xi}
\def\t{\theta}

\def\k{\kappa}

\def\a{\alpha}
\def\b{\beta}
\def\g{\gamma}
\def\d{\delta}
\def\l{\lambda}
\def\f{\frac}

\def\r{\rho}

\def\s{\sigma}

\def\di{\displaystyle}
\def\i{\infty}

\newtheorem{theorem}{Theorem}[section]
\newtheorem{lemma}[theorem]{Lemma}

\newtheorem{proposition}{Proposition}[section]
\newtheorem{remark}{Remark}

\renewcommand{\theequation}{\thesection.\arabic{equation}}

\begin{document}
	
	\title{\bf  Existence of Large Boundary Layer Solutions  to Inflow Problem of 1D Full 
		Compressible Navier-Stokes Equations} \vskip 0.5cm
	\author{Yi Wang\thanks{State Key Laboratory of Mathematical Sciences and Institute of Applied Mathematics, Academy of Mathematics and Systems Science, Chinese Academy of Sciences, Beijing 100190, P. R. China and School of Mathematical Sciences, University of Chinese Academy of Sciences, Beijing 100049, P. R. China ({\tt wangyi@amss.ac.cn}). The work of Yi Wang is partially supported by NSFC grants (Grant No.s 12171459, 12288201, 12090014, 12421001) and CAS Project for Young Scientists in Basic Research, Grant No. YSBR-031.}~, \quad Yong-Fu Yang\thanks{School of Mathematics, Hohai University, Nanjing 211100, P. R. China ({\tt yyang@hhu.edu.cn}).}~,\quad and Qiuyang Yu\thanks{Institute of Applied Mathematics, Academy of Mathematics and Systems Science, Chinese Academy of Sciences, Beijing 100190, P. R. China and School of Mathematical Sciences, University of Chinese Academy of Sciences, Beijing 100049, P. R. China ({\tt yuqiuyang@amss.ac.cn}).}
	}

	\date{Dedicated to Professor Boling Guo for His 90th Birthday}
	\maketitle
	
	\begin{abstract}
		We present the existence/non-existence criteria for large-amplitude boundary layer solutions to the inflow problem of the one-dimensional (1D) full compressible Navier--Stokes equations on a half line $\mathbb{R}_+$. Instead of the classical center manifold approach for the existence of small-amplitude boundary layer solutions in the previous results, the delicate global phase plane analysis, based on the qualitative theory of ODEs, is utilized to obtain the sufficient and necessary conditions for the existence/non-existence of large boundary layer solutions to the half-space inflow problem when the right end state belongs to the supersonic, transonic, and subsonic regions, respectively, which completely answers the existence/non-existence of boundary layer solutions to the half-space inflow problem of 1D full compressible Navier--Stokes equations.

		\

	\end{abstract}

	\section{Introduction}
	\renewcommand{\theequation}{\arabic{section}.\arabic{equation}}
	\setcounter{equation}{0}
	
	In this paper, we consider the full (or non-isentropic) compressible Navier--Stokes(--Fourier) equations on a half line $\mathbb{R}_+:=(0,+\infty)$:
	\begin{equation}\label{(1.1)}
		\left\{
		\begin{array}{ll}
			\di \rho_{t}+(\rho u)_x=0, & t\in \mathbb{R}_+,~~ x\in \mathbb{R}_+, \\[2mm]
			\di (\rho u)_t+(\rho u^2 +p)_x=\mu u_{xx},  &\\[2mm]
			\di \left[ \rho\left( e+\f{u^2}{2}\right) \right]_t+\left[ \rho
			u\left( \t+\f{u^2}{2}\right) +pu\right] _x=\k\t_{xx}+\mu(uu_x)_x, &
		\end{array}
		\right.
	\end{equation}
	where $\rho=\rho(t,x)>0$ is the fluid density, $u=u(t,x)$ is the velocity, and $\t=\t(t,x)$ is the absolute temperature,  respectively, and $p=p(\rho,\t)$ is the pressure, $e=e(\rho,\t)$ is the internal energy, $\mu>0$ is the viscosity coefficient, and $\k>0$ is the heat conductivity coefficient. Here we consider the classical case that both $\mu$ and $\kappa$ are positive constants. For the ideal and polytropic gas, one has the state equation
	\begin{equation}\label{state}
		p=R\rho\t=A\rho^\gamma\exp \left( \frac{\gamma-1}{R}s\right) ,~~~e=\frac{R\t}{\gamma-1}+\mbox{const},
	\end{equation}
	where $s=s(t,x)$ is the entropy, $\gamma>1$ is the adiabatic exponent, and $A,R>0$ are gas constants.
	
	We investigate the initial and boundary value problem \eqref{(1.1)} with the initial values
	\begin{equation}\label{initial}
		\begin{array}{l}
			\di (\rho,u,\t)(t=0,x)=(\rho_0,u_0,\t_0)(x)\to(\rho_+,u_+,\t_+),~~{\rm as}\ \ x\to+\i,\\[2mm]
		\end{array}
	\end{equation}
	where $\di \inf_{x\in\mathbb R_+}(\rho_0,\t_0)(x)>0$ and $\rho_+>0$, $\t_+>0$ and $u_+$ are given constants, and one of the following three boundary conditions \cite{Matsumura-2001}
	\begin{equation}\label{boundarycases}
		\begin{array}{l}
			\di u(t,x=0)=u_->0, ~~~(\rho,\t)(t,x=0)=(\rho_-,\t_-),\quad\mbox{inflow problem},\\[2mm]
			\di u(t,x=0)=u_-=0, ~~~\t(t,x=0)=\t_-,~~~~~~~~\!~~~~~~~~\mbox{impermeable wall problem},\\[2mm]
			\di u(t,x=0)=u_-<0, ~~~\t(t,x=0)=\t_-,~~~~~~~~\!~~~~~~~~\mbox{outflow problem},\\[2mm]
		\end{array}
	\end{equation}
	where $\rho_->0$, $\t_->0$, and $u_-$ are prescribed constants. Note that in \eqref{boundarycases}, we only concern the simple case that the boundary values are constants for all $t>0$. Furthermore, the compatible conditions should be satisfied for the initial values $\eqref{initial}$ and the boundary values $\eqref{boundarycases}$ at the origin point $(0,0)$. For the impermeable wall problem and outflow problem, that is, when the boundary velocity $u_-\leq 0$, the
	boundary density $\rho(t, 0)$ can not be imposed, while for the inflow problem with $u_->0$, the
	boundary density $\rho(t, 0)=\rho_->0$ must be given due to the well-posedness of the hyperbolic equation $\eqref{(1.1)}_1$.

	The large-time behavior of the solutions to the Cauchy problem of \eqref{(1.1)} on the whole line $\mathbb{R}$ with two different far-field states $(\rho_\pm, u_\pm, \theta_\pm)$ is expected to be determined by the Riemann solutions to the corresponding inviscid Euler equations
	\begin{equation}\label{Eulerequation}
		\left\{
		\begin{array}{l}
			\di \rho_{t}+(\rho u)_x=0,\\[2mm]
			\di (\rho u)_t+(\rho u^2 +p)_x=0,  \\[3mm]
			\di \left[ \rho\left( e+\f{u^2}{2}\right) \right]_t+\left[ \rho
			u\left( \t+\f{u^2}{2}\right) +pu\right] _x=0,
		\end{array}
		\right.
	\end{equation}
	with Riemann initial  data
	\begin{equation}\label{ri}
		(\rho,u,\t)(t=0,x)=\left\{
		\begin{array}{l}
			\di (\rho_-,u_-,\t_-),~~x<0,\\[2mm]
			\di (\rho_+,u_+,\t_+),~~x>0.
		\end{array}
		\right.
	\end{equation}
	The above Riemann solutions generically contain two nonlinear waves, i.e., shock and rarefaction waves in the genuinely nonlinear characteristic fields, and one linear wave, called contact discontinuity in the linearly degenerate characteristic field \cite{Smoller}. Generic Riemann solutions contain all of three basic wave patterns which are self-similar and dilation invariant, and govern not only the large-time behavior of solutions to compressible Navier--Stokes equations \eqref{(1.1)}, but also the local, global and large-time behavior of solutions to the compressible Euler system \eqref{Eulerequation} at least in small BV (bounded variation) regime.  The mathematical verification of the large-time behavior of the solutions to the Cauchy problem of \eqref{(1.1)} (and its isentropic simplification) towards the generic Riemann solutions, or in other words, the time-asymptotic stability of Riemann solutions to the Cauchy problem of compressible Navier--Stokes equations \eqref{(1.1)} are extensively investigated  since the pioneer work of Il'in--Oleinik \cite{Il'in-Oleinik} on Burgers equation, and much progress are made in this field by the direct or weighted energy methods, spectral methods, point-wise Green function methods, $L^1$-stability and even the combined methods mentioned above, see \cite{Goodman, MN-1985, MN-1986, Liu, LX, Szepessy-Xin, FD, MZ, Nishihara-Yang-Zhao, Huang-Xin-Yang, Huang-Li-Matsumura, LZ,HKK} and references therein. Very recently, Kang--Vasseur--Wang proved the generic Riemann solutions (containing viscous shock wave, rarefaction wave, and even viscous contact wave) to both barotropic compressible Navier--Stokes equations \cite{Kang-Vasseur-Wang-2023} and full compressible Navier--Stokes(--Fourier) equations \eqref{(1.1)} \cite{Kang-Vasseur-Wang-2024} with the help of $a$-contraction method for the stability of viscous shock wave with time-dependent shift \cite{KV}.

	For the initial and boundary value half-space problem \eqref{(1.1)}, due to the boundary effect, the boundary layer solution (stationary wave) may appear for the large-time behavior of the solution to \eqref{(1.1)}--\eqref{boundarycases}, besides the three basic wave patterns mentioned above. Matsumura \cite{Matsumura-2001} proposes a criterion when the boundary layer solution forms for the large-time behavior of the solution to barotropic Navier--Stokes equations. The argument in \cite{Matsumura-2001} for the barotropic Navier–Stokes equations can also be applied to the boundary layer solution to the half-space problem of full Navier--Stokes(--Fourier) equations \eqref{(1.1)}, see \cite{Huang-Li-Shi}. Consider the Riemann problem to the Euler equations \eqref{Eulerequation}, where the initial right-end state of the Riemann data in \eqref{ri} is given by the far field state $(\rho_+,u_+,\t_+)$ in \eqref{initial}, and the left-end state $(\rho_-,u_-,\theta_-)$ is given by all possible states which are consistent with the boundary conditions \eqref{boundarycases} on ${x=0}$. When the left-end state to the above Riemann problem is uniquely determined and consistent with the boundary conditions \eqref{boundarycases} on the boundary $x=0$, the boundary layer solution does not appear, that is, there is no need for the boundary layer. On the other hand, if the above Riemann solution on the boundary is not consistent with the boundary conditions \eqref{boundarycases}  for any admissible left-end states, we expect that a boundary layer solution comes up to compensate the gap. According to this criteria, for the impermeable wall problem, the boundary layer solution does not appear for the large-time behavior of the solution to both barotropic and full compressible Navier–Stokes equations, and for inflow/outflow problem, the boundary layer solution may occur. When the boundary layer solution appears, it can be constructed through the stationary solution to Navier--Stokes equations \eqref{(1.1)}. The existence and stability of the boundary layer solution to the inflow/outflow problem for the barotropic or full Navier--Stokes(--Fourier) equations have been studied extensively by many authors, see \cite{Matsumura-Nishihara-2001, KNZ, Huang-Matsumura-Shi, KZ-08, KZ-09, Huang-Qin, Qin-Wang-2009, Huang-Li-Shi, KNNZ, Qin-Wang-2011} and the references therein.
	
	As mentioned above, the boundary layer solution is constructed through the stationary solution to the half-space problem of the Navier--Stokes equations \eqref{(1.1)}. In the barotropic regime when the pressure $p$ depends only on the density function,  the energy equation $\eqref{(1.1)}_3$ can be dropped and then the mass equation and the momentum equation are combined into a closed system. Therefore, the stationary solution to barotropic Navier--Stokes equations $\eqref{(1.1)}_1$--$\eqref{(1.1)}_2$, after integrating over $[x, +\infty)$, can be simplified to an algebraic equation plus a first-order ODE equation on the half line $\mathbb{R}_+$ with the far-field condition and the boundary condition on $x=0$, which can be solved for arbitrarily large amplitude \cite{Matsumura-Nishihara-2001}.  While for full compressible Navier--Stokes(--Fourier) equations \eqref{(1.1)}, its stationary solution equation, after integrating over $[x, +\infty)$, can be transformed into an algebraic equation plus a first-order ODE system including the momentum and energy equations on the half line $\mathbb{R}_+$ with the far-field conditions and the boundary conditions on $x=0$. To solve this first-order ODE system with the algebraic constraint, center manifold  approach can be used for the existence of small-amplitude boundary layer solution, see \cite{Huang-Li-Shi, Qin-Wang-2009, Qin-Wang-2011} in Lagrangian coordinate and \cite{Nakamura-Nishibata} Eulerian coordinate, for details. Note that the small-amplitude conditions for the boundary layer solution are essential for the center manifold  approach  (see \cite{Bressan-2007,Carr}).
	
	In this paper we are concerned with the existence/non-existence criteria for large-amplitude boundary layer solution to the inflow problem  \eqref{(1.1)}--\eqref{initial}, $\eqref{boundarycases}_1$. Thanks to the fact that the stationary boundary layer solutions to  \eqref{(1.1)}--\eqref{initial}, $\eqref{boundarycases}_1$, after integrating over $[x, +\infty)$, satisfy the first-order planar autonomous ODE system, then by phase plane analysis as for the existence of large-amplitude shock profile to full compressible Navier--Stokes equations in Gilbarg \cite{Gilbarg}, we can present the sufficient and necessary conditions for the existence/non-existence of large-amplitude boundary layer solutions to the half-space inflow problem when the right end state belongs to the supersonic, transonic, and subsonic regions, respectively, which completely answer the existence/non-existence of boundary layer solutions to the inflow problem  \eqref{(1.1)}--\eqref{initial}, $\eqref{boundarycases}_1$.

	The rest of the paper will be organized as follows. In Section 2, we demonstrate the main result on the sufficient and necessary conditions for the existence/non-existence of large-amplitude boundary layer solutions to the inflow problem \eqref{(1.1)}-\eqref{initial}, $\eqref{boundarycases}_1$, and its detailed proof is given in Section 3. At last, we will give some comments on the time-asymptotic stability of these large-amplitude boundary layer solutions, and the corresponding results on the the existence of large-amplitude boundary layer solutions to the outflow problem \eqref{(1.1)}--\eqref{initial}, $\eqref{boundarycases}_3$. 
	
	\section{Main result}
	\setcounter{equation}{0}
	
	The main result will be shown in this section. Since the problem considered here is one-dimensional half-space inflow problem, we can rewrite the system \eqref{(1.1)} in Lagrangian coordinates by the coordinate transformation
	\begin{equation}\label{Lag}
		(t,x)\mapsto \left( t,\int_{(0,0)}^{(t,x)} \rho(\tau,y)dy-(\r
		u_1)(\tau,y)d\tau\right),
	\end{equation}
	where $\int_A^B fdy+gd\tau$ represents a line integration from the point $A$ to the point $B$ on
	${\mathbb R}_+\times {\mathbb R}_+$. Here, the line integration in \eqref{Lag} is independent of the integration path and then is unique due to the mass conservation laws.
	
	We will still denote the Lagrangian coordinates by $(t,x)$ for the simplicity of notations. Then the system \eqref{(1.1)} can be transformed into the following moving boundary problem in Lagrangian coordinates \cite{Matsumura-Nishihara-2001}:
	\begin{equation}\label{Lag-equation}
		\left\{
		\begin{array}{llll}
			\di v_t-u_{x}=0,\\[2mm]
			\di u_{t}+p_x=\mu\left( \f {u_x}v\right) _{x},~~~~~~~~~~~~~~~~~~~~~~~~~~~~~~~~~~~~~~~~~x>\s_-t,~~~t>0,\\[4mm]
			\di\left( e+\f{u^{2}}{2}\right) _{t}+
			(pu)_{x}=\k\left( \f{\t_x}{v}\right) _x+\mu\left( \f{uu_x}{v}\right) _x,\\[4mm]
			\di (v,u,\t)(t=0,x)=(v_0,u_0,\t_0)(x)\to(v_+,u_+,\t_+),~~\!~~\mbox{as}~ x\to+\i,\\[2mm]
			\di(v,u,\t)(t,x=\s_-t)=(v_-,u_-,\t_-),
		\end{array}
		\right. 
	\end{equation}
	where $v(t,x)=\frac{1}{\rho(t,x)}$ is the specific volume of the gas and the boundary moving speed $\s_-=-\frac{u_-}{v_-}<0$. In order to fix the moving boundary $x-\s_-t$, we introduce a new variable $\xi=x-\s_-t$. Then the moving boundary problem \eqref{Lag-equation} becomes a new half-space problem in Lagrangian coordinate:
	\begin{equation}\label{xi-equation}
		\left\{
		\begin{array}{llll}
			\di v_t-\s_-v_\xi-u_\xi=0,\\[2mm]
			\di u_{t}-\s_-u_\xi+p_\xi=\mu\left( \f {u_\xi}v\right) _{\xi},~~~~~~~~~~~~~~~~~~~~~~~~~~~~~~~~~~~~~~~~~~~~~~~~~~~~~~~\xi>0,~~~t>0,\\[4mm]
			\di\left(\frac{R\t}{\gamma-1}+\f{u^{2}}{2}\right) _{t}-\s_-\left(\frac{R\t}{\gamma-1}+\f{u^{2}}{2}\right) _\xi+
			(pu)_{\xi}=\k\left( \f{\t_\xi}{v}\right) _\xi+\mu\left( \f{uu_\xi}{v}\right) _\xi,\\[4mm]
			\di (v,u,\t)(t=0,\xi)=(v_0,u_0,\t_0)(\xi)\to(v_+,u_+,\t_+),~~~~~~~~~~~~~~~~~~~~~~~~~\!~~\mbox{as}~ \xi\to+\i,\\[2mm]
			\di(v,u,\t)(t,\xi=0)=(v_-,u_-,\t_-).
		\end{array}
		\right. 
	\end{equation}
	
	According to \cite{Matsumura-2001}, 
	the stationary boundary layer solution $(V,U,\Theta)(\xi)$ to \eqref{xi-equation} satisfies the following ODE system on $\mathbb{R}_+$:
	\begin{equation}\label{BL-equation}
		\left\{
		\begin{array}{llll}
			\di -\s_-V_\xi-U_\xi=0,\\[2mm]
			\di -\s_-U_\xi+P_\xi=\mu\left( \f {U_\xi}V\right) _{\xi},~~~~~~~~~~~~~~~~~~~~~~~~~~~~~~~~~~~~~~~~~~~~~~~~~\xi>0,\\[4mm]
			\di-\s_-\left(\frac{R\Theta}{\gamma-1}+\f{U^{2}}{2}\right) _\xi+
			(PU)_{\xi}=\k\left( \f{\Theta_\xi}{V}\right) _\xi+\mu\left( \f{UU_\xi}{V}\right) _\xi,\\[4mm]
			\di (V,U,\Theta)(0)=(v_-,u_-,\t_-),\\[2mm]
			\di(V,U,\Theta)(+\i)=(v_+,u_+,\t_+),
		\end{array}
		\right. 
	\end{equation}
	where $P=p(V,\Theta)=\frac{R\Theta}{V}$. The aim of the present paper is to present the existence/non-existence criteria for large-amplitude boundary layer solutions $(V,U,\Theta)(\xi)$ to ODE system \eqref{BL-equation} on the half line $\mathbb{R}_+$.	
	
	Define the local  sound speed $c(v,\t)$ and  Mach number $M(v,u,\t)$ as
	\begin{equation*}
		c(v,\t):=v\sqrt{\frac{\gamma p}{v}}=\sqrt{R\gamma\t},
	\end{equation*}
	and
	\begin{equation*}
		M(v,u,\t):=\frac{|u|}{c(v,\t)}=\frac{|u|}{\sqrt{R\gamma \t}}.
	\end{equation*}
	Our main result can be stated as follows.
	\begin{theorem}\label{maintheorem}
		For the inflow problem \eqref{Lag-equation} with $v_\pm>0$, $u_->0$, $\t_\pm>0$. 
		\begin{itemize}
			
			\item If $u_+\leq0$, then there is no boundary layer solution to \eqref{BL-equation}. 
			
			\item If $u_+>0$, then the existence/nonexistence of the boundary layer solutions to \eqref{BL-equation} can be divided into the following three cases according to the Mach number $M_+:=\frac{|u_+|}{\sqrt{R\gamma \t_+}}$:
			
			\underline{Case I} (Supersonic case): $M_+>1$. Then there is no solution to \eqref{BL-equation}.
			
			\underline{Case II} (Transonic case): $M_+=1$. Then there exists a unique curve $\Sigma$ (see Figure \ref{Fig-1} below) such that the following holds.
			
			\begin{itemize}
				\item [(1).] $\Sigma$ is monotone in $(u,\t)$-space with $0<u<u_+$ and $0<\theta_+<\theta$.
				
				\item [(2).] $\Sigma$ is tangent to the straight line
				\begin{equation*}
					(\gamma-1)u_+(u-u_+)+R\gamma(\t-\theta_+)=0
				\end{equation*}
				in $(u,\t)$-space at the equilibrium point $(u_+,\t_+)$.
				
				\item[(3).] The boundary layer solution $(V,U,\Theta)(\xi)$ to \eqref{BL-equation} exists if and only if $(u_-,\t_-)\in \Sigma$ and $\frac{u_-}{v_-}=\frac{u_+}{v_+}$. Moreover, the solution $(V,U,\Theta)(\xi)$ is monotone $(V'>0,U'>0,\Theta'<0)$ and satisfies
				\begin{equation}\label{trans-decay}
					\left|\frac{d^n}{d\x^n}(V-v_+,U-u_+,\Theta-\t_+) \right| \leq  \frac{C}{(1+\xi)^{n+1}},\quad \forall \xi>0, \ n=0,1,2,\cdots,
				\end{equation}
				for some positive constant $C$.
			\end{itemize}
			
			\underline{Case III} (Subsonic case): $0<M_+<1$. Then there exists a unique piecewise smooth curve consisting of two segments $\Gamma_1$ and $\Gamma_2$ (see Figure \ref{Fig-2} below) such that the following holds.
			
			\begin{itemize}
				\item [(1).]  Both $\Gamma_1$ and $\Gamma_2$ are monotone in $(u,\t)$-space and $\Gamma_1$ is defined on $0<u<u_+$ and $\Gamma_2$ is defined on $\max \left\lbrace 0, \a_2\t_+\right\rbrace <\t<\t_+$ with $\a_2:=1-\frac{2(1-M_+^2)(M_+^2\g+1)(\g-1)}{M_+^2(\g+1)^2 }$.
				
				\item [(2).] Both $\Gamma_1$ and $\Gamma_2$ are tangent to the straight line
				\begin{equation*}
					u_+^2(u-u_+)+M_+^2\gamma\kappa\left(\frac{Ru_+}{\k(\gamma-1)} -\l_2\right) (\t-\theta_+)=0
				\end{equation*}
				in $(u,\t)$-space at the equilibrium point $(u_+,\t_+)$, where $\l_2$ is the negative eigenvalue of the matrix $A$ (to be defined in \eqref{A}).
				
				\item[(3).] The boundary layer solution $(V,U,\Theta)(\xi)$ to \eqref{BL-equation} exists if and only if $(u_-,\t_-)\in \Gamma_1\cup\Gamma_2$ and $\frac{u_-}{v_-}=\frac{u_+}{v_+}$. The solution $(V,U,\Theta)(\xi)$ is monotone with $V'>0$, $U'>0$, $\Theta'<0$ if $(u_-,\t_-)\in \Gamma_1$,  and  $V'<0$, $U'<0$, $\Theta'>0$  if $(u_-,\t_-)\in \Gamma_2$. Moreover, the solution $(V,U,\Theta)(\xi)$ satisfies
				\begin{equation}\label{subs-decay}
					\left|\frac{d^n}{d\xi^n}(V-v_+,U-u_+,\Theta-\t_+) \right| \leq Ce^{-c\xi},\quad \forall \xi>0, \ n=0,1,2,\cdots,
				\end{equation}
				for some positive constants $C$ and $c$.
			\end{itemize}
		\end{itemize}
	\end{theorem}
	
	\begin{remark}
		Theorem \ref{maintheorem} does not require any smallness conditions on the amplitude of the boundary layer solution, which completely answers the existence/non-existence of boundary layer solutions \eqref{BL-equation} to the inflow problem  \eqref{(1.1)}--\eqref{initial}, $\eqref{boundarycases}_1$ or \eqref{Lag-equation}.
	\end{remark}
	\begin{figure}[htbp]
		\centering
		\includegraphics[width=0.5\textwidth]{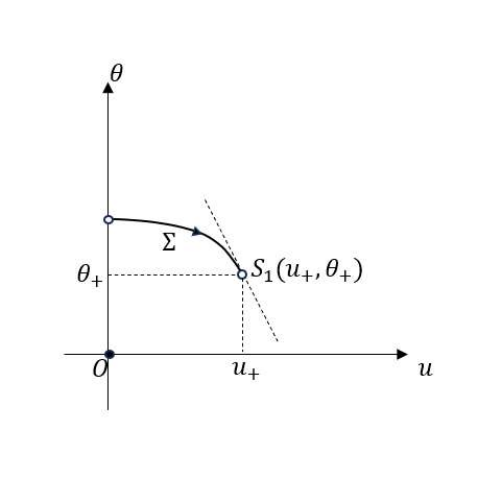}
		\caption{Transonic case.}
		\label{Fig-1}
	\end{figure}
	\begin{figure}[htbp]
		\centering
		\begin{minipage}[t]{0.45\textwidth}
			\centering
			\includegraphics[width=\textwidth]{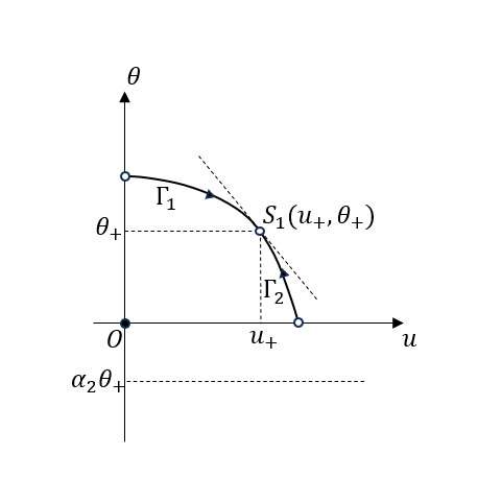}
			\label{fig:image1}
		\end{minipage}
		\hfill
		\begin{minipage}[t]{0.45\textwidth}
			\centering
			\includegraphics[width=\textwidth]{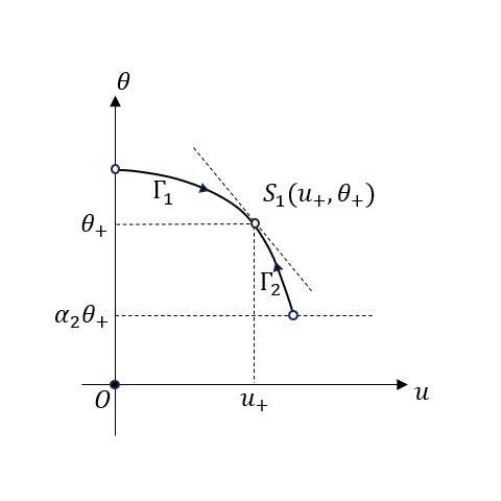}
			\label{fig:image2}
		\end{minipage}
		\caption{Subsonic case.}
		\label{Fig-2}
	\end{figure}
	
	\

	\section{Proof of main result}
	\setcounter{equation}{0}
	
	In this section we prove Theorem \ref{maintheorem}. Integrating the ODE system \eqref{BL-equation} over $[\xi,+\i)$ yields the following algebraic equation plus a first-order planar autonomous system
	\begin{equation}\label{integration}
		\left\{
		\begin{array}{llll}
			\di -\s_-(V-v_+)-(U-u_+)=0,\\[2mm]
			\di \mu \f {U'}V=-\s_-(U-u_+)+R\left( \frac{\Theta}{V}-\frac{\theta_+}{v_+}\right),\\[4mm]
			\di \k\f {\Theta'}V=-\s_-\frac{R}{\gamma-1}(\Theta-\t_+)+p_+(U-u_+)+\frac{\s_-}{2}(U-u_+)^2,\\[4mm]
			\di (V,U,\Theta)(0)=(v_-,u_-,\t_-),\\[2mm]
			\di(V,U,\Theta)(+\i)=(v_+,u_+,\t_+),
		\end{array}
		\right. 
	\end{equation}
	where $^\prime =\frac{d}{d\x}$, $p_+=p(v_+,\t_+)=\frac{R\theta_+}{v_+}$. Set $\xi=0$ in the algebraic equation $\eqref{integration}_1$, we have
	$$\s_-=-\frac{u_-}{v_-}=-\frac{u_+}{v_+}=-\frac{U(\xi)}{V(\xi)}.$$
	Thus
	\begin{equation}
		u_+=\frac{u_-}{v_-}v_+>0,
	\end{equation}
	which is a necessary condition for the existence of the boundary layer solution  $(V,U,\Theta)(\xi)$. More precisely, we require that
	\begin{equation}\label{Upositive}
		U(\xi)=\frac{u_-}{v_-}V(\xi)>0,~~~~~~\forall \xi\in[0,+\i).
	\end{equation}
	Furthermore, to give the solution physical meaning, we require
	\begin{equation}\label{thermodynamicslaw}
		\Theta(\xi)>0,~~~~\xi\in[0,+\i),~~~~\mbox{the third law of thermodynamics}.
	\end{equation}
	
	Plugging $V(\xi)=\frac{v_-}{u_-}U(\xi)$ into $\eqref{integration}_2$, $\eqref{integration}_3$, and then linearizing the right-hand side of $\eqref{integration}_2$ and $\eqref{integration}_3$ around $(u_+,\t_+)$, one can get a closed first-order planar autonomous system
	\begin{equation}\label{keyequation}
		\left\{
		\begin{array}{ll}
			\begin{pmatrix}
				\di U\\
				\di \Theta
			\end{pmatrix}'
			=
			A
			\begin{pmatrix}
				\di U-u_+\\
				\di \Theta-\t_+
			\end{pmatrix}
			+
			\begin{pmatrix}
				\di F_1(U,\Theta)\\
				\di F_2(U,\Theta)
			\end{pmatrix},\\[5mm]
			(U,\Theta)(0)=(u_-,\t_-),\\[3mm]
			(U,\Theta)(+\i)=(u_+,\t_+),
		\end{array}
		\right.
	\end{equation}
	with the matrix
	\begin{equation}\label{A}
		A:=
		\begin{pmatrix}
			\di \frac{(M_+^2\gamma-1)u_+}{M_+^2\gamma\mu} & ~\di\frac{R}\mu\\[5mm]
			\di \frac{u_+^2}{M_+^2\g\k} & ~\di\frac{Ru_+}{\k(\g-1)}
		\end{pmatrix},
	\end{equation}
	and the nonlinear terms
	\begin{equation}
		F_1(U,\Theta):=\f1\mu(U-u_+)^2,
	\end{equation}
	\begin{equation}
		\begin{array}{ll}
			\di F_2(U,\Theta):=\left(\frac{R\theta_+}{\k u_+}-\frac{u_+}{2\k} \right)(U-u_+)^2+\frac{R}{\k(\g-1)}(U-u_+)(\Theta-\theta_+)-\f1{2\k}(U-u_+)^3.
		\end{array}
	\end{equation}
	
	To prove Theorem \ref{maintheorem}, it is equivalent to consider the first-order planar autonomous system \eqref{keyequation}. We first give a proposition.
	
	\begin{proposition}\label{s-n-c}
		The solution to \eqref{BL-equation} for the inflow problem $(u_->0)$ exists if and only if the following both conditions hold.
		
		\
		
		(A). $\di\frac{u_-}{v_-}=\frac{u_+}{v_+}$.
		
		(B). The solution to \eqref{keyequation} exists with
		\begin{equation*}
			U(\xi)>0,~~\Theta(\xi)>0,~~~~\forall\xi\in[0,+\i).
		\end{equation*}
	\end{proposition}
	
	Since $v_\pm>0$ and the inflow boundary condition $u_->0$, it is obvious that $u_+\leq0$ violates the above condition (A), and then the solution to \eqref{BL-equation} does not exist if  $u_+\leq0$. Therefore,  we only focus on the case  $u_+>0$ to the system \eqref{keyequation} in the sequel.
	
	Direct computations show that the ODE system $\eqref{keyequation}_1$ has three equilibrium points 
	$$O(0,0),~~S_1(u_+,\t_+),~~S_2\left(\a_1u_+,\a_2\t_+ \right),$$
	where 
	\begin{equation}\label{a1}
		\a_1:=\frac{M_+^2\g-M_+^2+2}{M_+^2(\g+1)},
	\end{equation}
	and 
	\begin{equation}\label{a2}
		\a_2:=1-\frac{2(1-M_+^2)(M_+^2\g+1)(\g-1)}{M_+^2(\g+1)^2 }.
	\end{equation}
	It is worth noting that for the transonic case when $M_+=1$, the last two equilibrium points $S_1$ and $S_2$ coincide. Since $(U,\Theta)(\xi)$ tends to the equilibrium point $S_1(u_+,\t_+)$ as $\xi\to+\i$, we need focus on the classification of the equilibrium point $S_1(u_+,\t_+)$. To this end, we find that the determinant of the matrix $A$ is 
	\begin{equation}\label{determinant}
		\det A =\left| \begin{matrix}
			\di \frac{(M_+^2\gamma-1)u_+}{M_+^2\gamma\mu} & ~\di\frac{R}\mu\\[4mm]
			\di \frac{u_+^2}{M_+^2\g\k} & ~\di\frac{Ru_+}{\k(\g-1)}
		\end{matrix}\right| 
		=\frac{Ru_+^2(M_+^2-1)}{M_+^2\mu\kappa (\g-1)}.
	\end{equation}

	Then the existence/non-existence of solutions to \eqref{BL-equation} is divided into the following three cases according to the sign of the determinant of $A$, or equivalently, the value of the far-field Mach number $M_+$.
	
	\
	
	\underline{Case I} (Supersonic case): $M_+>1$. Then $\det A>0$ and
	\begin{equation}\label{trace}
		\operatorname{tr} A=\frac{Ru_+}{\kappa(\g-1)}+\frac{(M_+^2-1)u_+}{M_+^2\mu}+\frac{(\g-1)u_+}{M_+^2\g\mu}=:\b_1+\b_2+\b_3>0.
	\end{equation}
	Then we can deduce that the characteristic polynomial of $A$
	$$\l^2-\l\operatorname{tr} A +\det A$$
	has two different positive eigenvalues since
	$$(\operatorname{tr} A)^2-4\det A=(\b_1+\b_2+\b_3)^2-4\b_1\b_2=(\b_1-\b_2)^2+\b_3^2+2\b_1\b_3+2\b_2\b_3>0.$$
	Thus the equilibrium point $(u_+,\t_+)$ is an unstable node of $\eqref{keyequation}_1$, which means that no integral curves of $\eqref{keyequation}_1$ can satisfy $(U,\Theta)(+\i)=(u_+,\t_+)$. In other words, the solution to \eqref{keyequation} does not exist.

	\

	\underline{Case II} (Transonic case): $M_+=1$. Then $\det A=0$. It is easy to calculate that the two eigenvalues and the corresponding eigenvectors of $A$ are
	\begin{equation*}
		\begin{matrix}
			\di \l_1=0, & ~~~\di e_1=\left(1,-\frac{(\g-1)u_+}{R\g} \right),\\[4mm]
			\di \l_2=\left( \frac{\g-1}{\g\mu}+\frac{R}{\k(\g-1)}\right) u_+>0, & ~~~\di  e_2=\left(1,\frac{\mu u_+}{\k(\g-1)} \right).
		\end{matrix}
	\end{equation*}
	This case is the most complicated one, as one of the eigenvalues of $A$ is zero and  the equilibrium point $S_1(u_+,\t_+)$ is degenerate. To diagonalize the system \eqref{keyequation}, we set
	\begin{equation*}
		P:=
		\begin{pmatrix}
			\di 1 & ~~~\di 1\\[4mm]
			\di -\frac{(\g-1)u_+}{R\g} & ~~~\di \frac{\mu u_+}{\k(\g-1)}
		\end{pmatrix}.
	\end{equation*}
	Note that the two columns of $P$ are two eigenvectors of $A$ respectively, and $P$ is invertible since $e_1$ and $e_2$ are linearly independent. Direct calculations show that
	\begin{equation*}
		P^{-1}=\frac{1}{\det P}
		\begin{pmatrix}
			\di \frac{\mu u_+}{\k(\g-1)} & ~\di -1\\[4mm]
			\di \frac{(\g-1)u_+}{R\g} & ~\di 1
		\end{pmatrix}
	\end{equation*}
	satisfying
	\begin{equation*}
		P^{-1}AP=
		\begin{pmatrix}
			\di 0 & \di 0\\
			\di 0 & \di \l_2
		\end{pmatrix}.
	\end{equation*}
	Therefore we can perform the change of variables
	\begin{equation*}
		W=
		\begin{pmatrix}
			\di W_1 \\
			\di W_2
		\end{pmatrix}
		:=P^{-1}\begin{pmatrix}
			\di U-u_+ \\
			\di \Theta-\t_+
		\end{pmatrix},
	\end{equation*}
	such that \eqref{keyequation} can be rewritten as
	\begin{equation}\label{W}
		\left\{
		\begin{array}{llll}
			\di W_{1\xi}=g_1(W_1,W_2),\\[4mm]
			\di W_{2\xi}=\l_2W_2+g_2(W_1,W_2),\\[5mm]
			\di \begin{pmatrix}
				\di W_1 \\
				\di W_2
			\end{pmatrix}
			(0)=P^{-1}\begin{pmatrix}
				\di u_--u_+ \\
				\di \t_--\t_+
			\end{pmatrix},\\[5mm]
			\di \begin{pmatrix}
				\di W_1 \\
				\di W_2
			\end{pmatrix}(+\i)=\begin{pmatrix}
				\di 0 \\
				\di 0
			\end{pmatrix},
		\end{array}
		\right. 
	\end{equation}
	where
	\begin{equation*}
		\begin{array}{ll}
			\di g_1(W_1,W_2)=\frac{1}{\det P}\Bigg[ \left( \frac{u_+}{\k(\g-1)}-\frac{(2-\g)u_+}{2\g\k}\right) (W_1+W_2)^2\\[5mm]
			\di\qquad-\frac{R}{\k(\g-1)}(W_1+W_2)\left(-\frac{(\g-1)u_+}{R\g}W_1+\frac{\mu u_+}{\k(\g-1)}W_2 \right) +\frac{1}{2\k}(W_1+W_2)^3\Bigg],
		\end{array}
	\end{equation*}
	\begin{equation*}
		\begin{array}{ll}
			\di g_2(W_1,W_2)=\frac{1}{\det P}\Bigg[ \left( \frac{(\g-1)u_+}{R\g\mu}+\frac{(2-\g)u_+}{2\g\k}\right) (W_1+W_2)^2\\[7mm]
			\di\qquad+\frac{R}{\k(\g-1)}(W_1+W_2)\left(-\frac{(\g-1)u_+}{R\g}W_1+\frac{\mu u_+}{\k(\g-1)}W_2 \right) -\frac{1}{2\k}(W_1+W_2)^3\Bigg].
		\end{array}
	\end{equation*}
	The qualitative theory of ODEs can characterize the solution behavior of \eqref{W} near the equilibrium point $(0,0)$, which can be summarized in the following lemma.
	
	\begin{lemma}\label{lemma3.1}
		(\cite{Zhang-Ding-Huang-Dong}) Consider the following planar autonomous system
		\begin{equation}\label{planar autonomous system}
			\left\{
			\begin{array}{llll}
				\di \frac{dx}{d\xi}=g_1(x,y),\\[4mm]
				\di \frac{dy}{d\xi}=\l y+g_2(x,y).
			\end{array}
			\right. 
		\end{equation}
		Assume that $(0,0)$ is an isolated equilibrium point of \eqref{planar autonomous system}, together with $\l>0$ and $g_1(x,y)$, $g_2(x,y)$ being analytic functions of order not less than 2 within $B(0,\d)$. Let $\phi=\phi(x)~(|x|<\d)$ be the implicit function determined by the equation $\l \phi+g_2(x,\phi)=0,~|x|<\d$. Denote 
		$$\psi(x):=g_1(x,\phi(x))=a_mx^m+o(x^m),$$
		where $a_m\neq 0$, $m\geq2$. Then the following properties hold.
		
		(1) Any integral curve of $\eqref{planar autonomous system}$ which approaches $(0,0)$ as $\xi\to+\i$ or $\xi\to-\i$ must be tangent to one of the coordinate axes, i.e. $x=0$ or $y=0$.
		
		(2) If $m$ is odd and $a_m>0$, then $(0,0)$ is an unstable node of \eqref{planar autonomous system}.
		
		(3) If $m$ is odd and $a_m<0$, then $(0,0)$ is a saddle point of \eqref{planar autonomous system}.
		
		(4) If $m$ is even, then $(0,0)$ is a saddle-node point of \eqref{planar autonomous system}. Specifically, if $a_m>0$, then there exists a unique solution to \eqref{planar autonomous system} satisfying $\di\lim_{\xi\to+\i}(x(\xi),y(\xi))=(0,0)$ and tangent to the negative half $x$-axis $\{(x,0)|x<0\}$ (see Figure 3(a)). If $a_m<0$, then there exists a unique solution to \eqref{planar autonomous system} satisfying $\di\lim_{\xi\to+\i}(x(\xi),y(\xi))=(0,0)$ and tangent to the positive half $x$-axis $\{(x,0)|x>0\}$ (see Figure 3(b)).
	\end{lemma}
	
	\begin{figure}[htbp]
		\centering
		\begin{subfigure}[t]{0.45\textwidth}
			\centering
			\includegraphics[width=\textwidth]{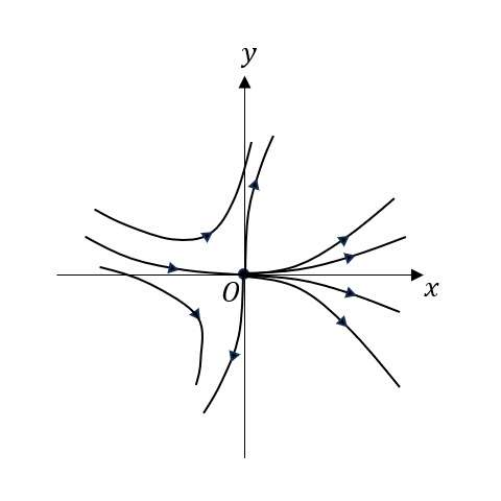}
			\caption{$a_m>0$.}
			\label{fig:sub3-1}
		\end{subfigure}
		\hfill
		\begin{subfigure}[t]{0.45\textwidth}
			\centering
			\includegraphics[width=\textwidth]{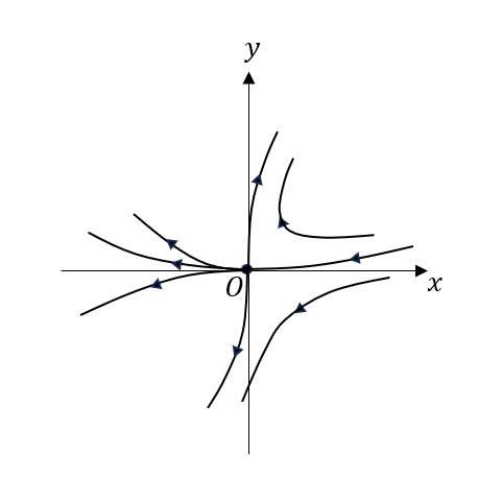}
			\caption{$a_m<0$.}
			\label{fig:sub3-2}
		\end{subfigure}
		\caption{Saddle-node point.}
		\label{Fig-3}
	\end{figure}

	We will use Lemma \ref{lemma3.1} to analyze \eqref{W}. Let $\phi=\phi(W_1)$ be the implicit function determined by the equation $\l_2 \phi+g_2(W_1,\phi)=0,~|W_1|<\d$. Then it is easy to check that $\phi(0)=\phi'(0)=0$. Direct calculations show that
	\begin{equation*}
		\begin{array}{ll}
			\di \psi(W_1):=g_1(W_1,\phi(W_1))=\frac{R\g(\g+1)}{2\left[R\g\mu+\k(\g-1)^2 \right] }W_1^2+o(W_1^2),
		\end{array}
	\end{equation*}
	which implies that, by Lemma  \ref{lemma3.1}  (4),  $(0,0)$ is a saddle-node point of \eqref{planar autonomous system} and there exists a unique integral curve $(W_1(\xi),W_2(\xi))$ of \eqref{planar autonomous system} satisfying $\di\lim_{\xi\to+\i}(W_1(\xi),W_2(\xi))=(0,0)$, which is tangent to the negative $W_1$-axis $\{(W_1,0)|W_1<0\}$. Back to the original $(u,\t)$-space, there exists a unique integral curve $(U(\xi),\Theta(\xi))$ of $\eqref{keyequation}_1$ satisfying $\di\lim_{\xi\to+\i}(U(\xi),\Theta(\xi))=(u_+,\t_+)$, which is tangent to the half straight line $$\tau:=\left\{(u,\t)|\di(\gamma-1)u_+(u-u_+)+R\gamma(\t-\theta_+)=0,~~u\leq u_+\right\}$$
	at the point $(u_+,\t_+)$. This curve $(U(\xi),\Theta(\xi))$ is denoted as $\Sigma$.
	
	The remaining task is to consider the direction of the trajectory $\Sigma$ in $(u,\t)$-space when $\xi$ decreases from $+\i$.  For this, consider two curves $\t=h_1(u)$ and $\t=h_2(u)$ in $(u,\t)$-space, where
	$$h_1(u):=-\frac{1}{R}(u-u_+)\left( u-\frac{1}{M_+^2\g}u_+\right)+\t_+,$$
	$$h_2(u):=\frac{\g-1}{2R}(u-u_+)\left( u-\frac{(M_+^2\g+2)u_+}{M_+^2\g} \right)+\t_+,$$
	which correspond to $\frac{dU}{d\xi}=0$ and $\frac{d\Theta}{d\xi}=0$, respectively. Set
	$$l_1:=\{(u,\t)|\t=h_1(u),~0\leq u\leq u_+\},$$
	\begin{equation}\label{RegionIboundary}
		l_2:=\{(u,\t)|\t=h_2(u),~0\leq u\leq u_+\},
	\end{equation}
	$$l_3:=\{(u,\t)|u=0\}.$$
	Let us denote the region enclosed by $l_1$, $l_2$ and $l_3$ as Region I (see Figure \ref{Fig-4} below).  Since $\frac{dU}{d\xi}=0$ and $\frac{d\Theta}{d\xi}\leq0$ on $l_1$, any integral curve $(U,\Theta)$ of $\eqref{keyequation}_1$ passing through the points on $l_1$ has a tangent vector pointing vertically downwards at those points (see Figure \ref{Fig-4}).	Similarly, any integral curve $(U,\Theta)$ of $\eqref{keyequation}_1$ passing through the points on $l_2$ has a tangent vector pointing horizontally to the right at those points. Since $\Sigma$ is tangent to the half straight line $\tau$ at the equlibrum point $S_1(u_+,\t_+)$ when $\xi \rightarrow +\infty$, we can assert from Figure \ref{Fig-4} that $\Sigma$ lies in Region I when $\xi$ is sufficiently large. Moreover, due to the fact that $\frac{d\Theta}{dU}<0$ in Region I, we can deduce that $\Sigma$ is monotone in $(u,\t)$-space. Therefore, when we extend $\Sigma$ as $\xi$ decreases from $+\i$, $\Sigma$ can not intersect with the boundaries $l_1$ or $l_2$,  or terminate inside Region I as $\xi\to-\i$ since there is no equilibrium point inside the Region I. As a consequence, $\Sigma$ has to intersect with the boundary $l_3$ at a unique point $Z_0\in l_3$. On the other hand, the necessary condition \eqref{Upositive} for the existence of the boundary layer solution is $U(\x)>0$, while $U(\xi)\equiv0$ on the whole line $l_3$ and $U(\xi)<0$ on the left of $l_3$. Thus, the curve $\Sigma$ should be restricted between the point $Z_0\in l_3$ and the equilibrium point $S_1(u_+,\t_+)$.
	
	In conclusion, the solution to \eqref{keyequation} exists if and only if the boundary data $(u_-,\t_-)\in \Sigma$. Together with the condition (A) in Proposition \ref{s-n-c}, we can prove Theorem \ref{maintheorem}  in transonic case.
	\begin{figure}[htbp]
		\centering
		\includegraphics[width=0.8\textwidth]{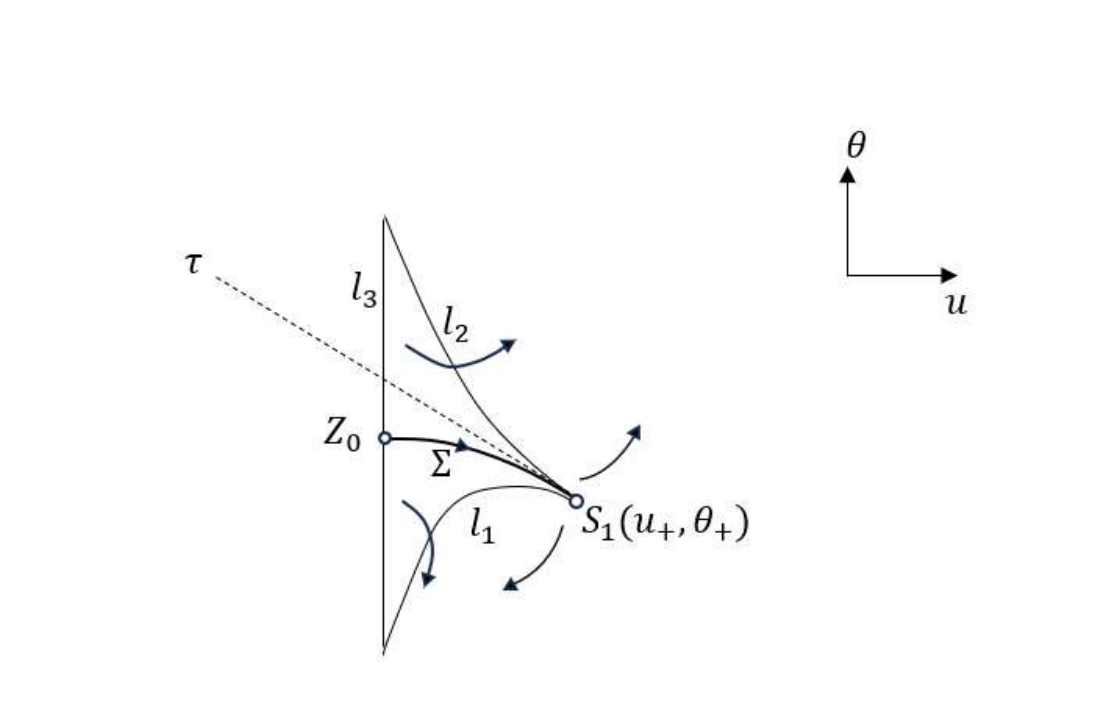}
		\caption{Region I in transonic case.}
		\label{Fig-4}
	\end{figure}
	
	\
	
\begin{figure}[htbp]
		\centering
		\begin{subfigure}[t]{0.5\textwidth}
			\centering
			\includegraphics[width=\textwidth]{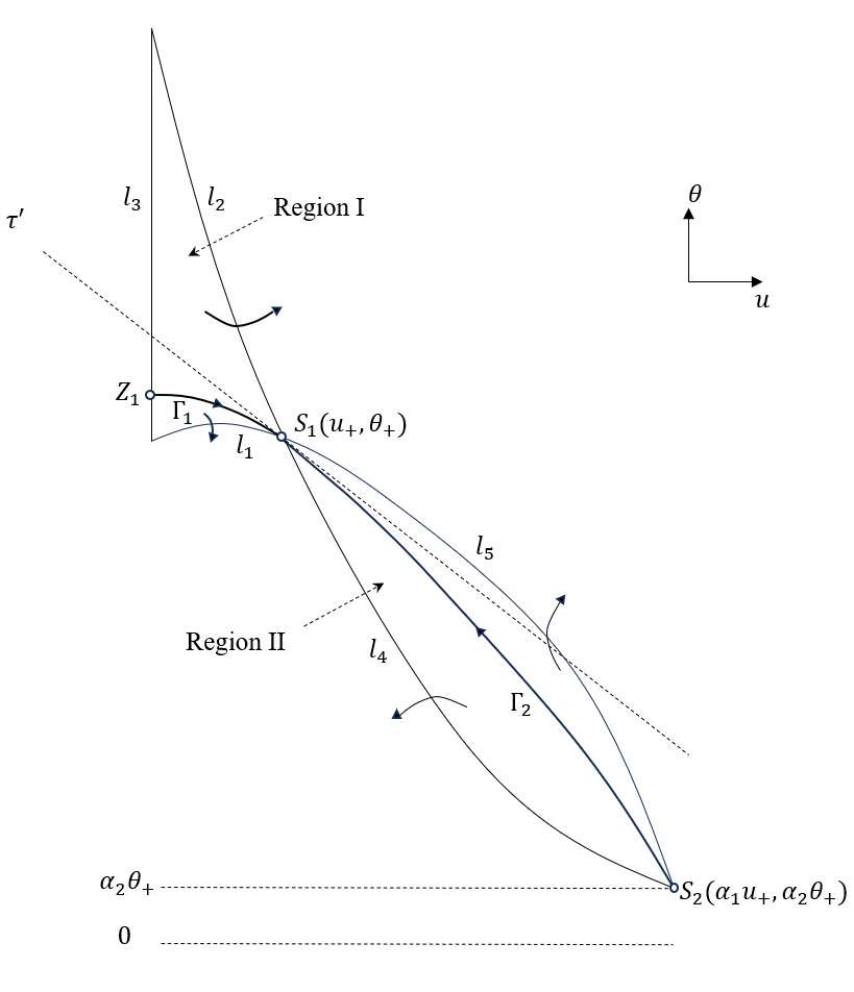}
			\caption{$\sqrt{\frac{\g-1}{2\g}}< M_+<1$.}
			\label{fig:sub5-1}
		\end{subfigure}
		\hfill
		\begin{subfigure}[t]{0.492\textwidth}
			\centering
			\includegraphics[width=\textwidth]{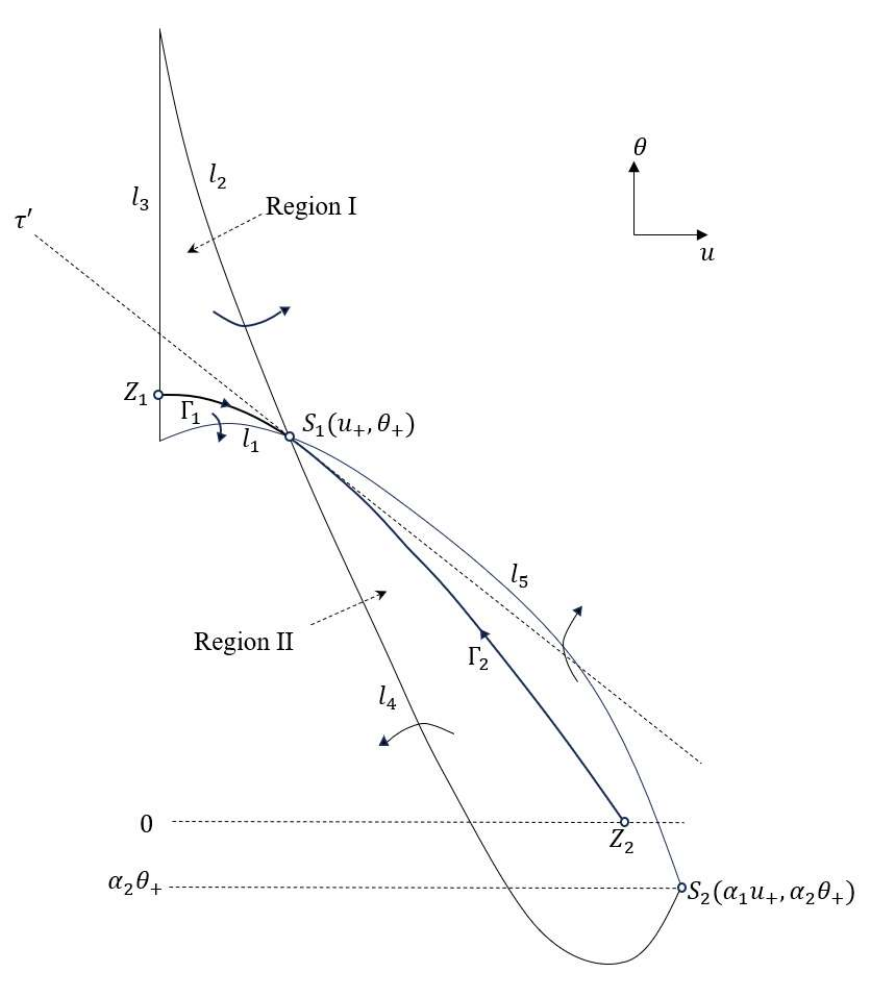}
			\caption{$0<M_+\leq\sqrt{\frac{\g-1}{2\g}}$.}
			\label{fig:sub5-2}
		\end{subfigure}
		\caption{Region I and Region II in subsonic case.}
		\label{Fig-5}
	\end{figure}

	\underline{Case III} (Subsonic case): $0<M_+<1$. Then $\det A<0$, and $A$ has two real eigenvalues $\l_1$ and $\l_2$ with the opposite signs.  For definiteness, assume that $\l_1>0$, $\l_2<0$. Therefore, the  equilibrium point $S_1(u_+,\t_+)$ is a saddle point of $\eqref{keyequation}_1$. According to the qualitative theory of ODEs, we know that only two integral curves, denoted by $\Gamma_1$ and $\Gamma_2$, tend to  the equilibrium point $S_1(u_+,\t_+)$ as $\xi\to+\i$ (see Figure \ref{Fig-5}). Both curves  $\Gamma_1$ and $\Gamma_2$ are tangent to the straight line $$\tau':=\left\{(u,\t)~\bigg|~\di u_+^2(u-u_+)+M_+^2\gamma\kappa\left(\frac{Ru_+}{\k(\gamma-1)} -\l_2\right) (\t-\theta_+)=0\right\}$$
	at the point $(u_+,\t_+)$.
	
	Using the similar arguments as in the transonic case (Case II), we investigate the tangent vectors of the integral curve of $\eqref{keyequation}_1$ in Region I and Region II in Figure \ref{Fig-5}. Let the boundary of Region I defined exactly same as in Case II in \eqref{RegionIboundary}, while Region II is enclosed by the two curves $l_4$ and $l_5$ defined by
	$$l_4:=\{(u,\t)|\t=h_2(u),~u_+\leq u\leq \a_1u_+\},$$
	$$l_5:=\{(u,\t)|\t=h_1(u),~u_+\leq u\leq \a_1u_+\},$$
	where $\a_1$ is given in \eqref{a1}. As $\xi\rightarrow +\infty$, both curves  $\Gamma_1$ and $\Gamma_2$ tend to the equilibrium point $S_1(u_+,\t_+)$ from Region I and II, respectively. Then we extend both curves  $\Gamma_1$ and $\Gamma_2$ as $\xi$ decreases from $+\i$.

	The extension of $\Gamma_1$ is similar to and even simpler than Case II.  Since $\frac{d\Theta}{dU}<0$ inside the whole Region I, the curve $\Gamma_1$ is strictly monotone in Region I and can not intersect with the boundaries $l_1$ and $l_2$.	Thus $\Gamma_1$ has to intersect with the boundary $l_3$ at a unique point $Z_1\in l_3$. Therefore, the curve $\Gamma_1$ should be restricted between the point $Z_1\in l_3$ and the equilibrium point $S_1(u_+,\t_+)$ as in Figure \ref{Fig-5}. Finally, if $(u_-,\t_-)\in\Gamma_1$ and $\frac{u_-}{v_-}=\frac{u_+}{v_+}$, then we can assert that the boundary layer solution to \eqref{BL-equation} exists. On the other hand, when we investigate the extension of the curve $\Gamma_2$ in Region II, we need first consider the relative positions of the equilibrium point  $S_2\left(\a_1u_+,\a_2\t_+ \right)$ and the horizontal line $\t=0$ to obey \eqref{thermodynamicslaw}.
	Thus we classify the considerations into two subcases.

	(a). $\sqrt{\frac{\g-1}{2\g}}< M_+<1$. In this case, the equilibrium point $S_2\left(\a_1u_+,\a_2\t_+ \right)$ is above the horizontal line $\t=0$ and both $l_4$ and $l_5$ are monotone in $(u,\t)$-space. Hence it is crucial to analyze the category of the equilibrium point $S_2$. Denote $(U',\Theta')=:(H_1(U,\Theta),H_2(U,\Theta))$ and expand the planar autonomous system $\eqref{keyequation}_1$ around the equilibrium point $S_2\left(\a_1u_+,\a_2\t_+ \right)$, one can get
	\begin{equation}\label{change}
		\begin{array}{ll}
			\begin{pmatrix}
				\di U\\[3mm]
				\di \Theta
			\end{pmatrix}'
			=
			B
			\begin{pmatrix}
				\di U-\a_1u_+\\[3mm]
				\di \Theta-\a_2\t_+
			\end{pmatrix}
			+
			\begin{pmatrix}
				\di G_1(U,\Theta)\\[3mm]
				\di G_2(U,\Theta)
			\end{pmatrix},
		\end{array}
	\end{equation}
	where the matrix $B=\left.\begin{pmatrix}
		\di \frac{\partial H_1}{\partial u} & ~\di\frac{\partial H_1}{\partial \theta}\\[4mm]
		\di \frac{\partial H_2}{\partial u} & ~\di\frac{\partial H_2}{\partial \theta}
	\end{pmatrix}\right|_{S_2} $ and $G_1$, $G_2$ are quadratic with respect to $(U-\a_1u_+, \Theta-\a_2\t_+)$  near the point $S_2\left(\a_1u_+,\a_2\t_+ \right)$. Since $\frac{\partial H_1}{\partial \theta}>0$ and $\frac{\partial H_2}{\partial u}>0$ at the point $S_2\left(\a_1u_+,\a_2\t_+ \right)$, it is easy to check that the matrix $B$ has two different positive eigenvalues. Therefore,  the equilibrium point $S_2\left(\a_1u_+,\a_2\t_+ \right)$ is an unstable node of $\eqref{keyequation}_1$ or \eqref{change}.
	
	Now we are at the stage to extend $\Gamma_2$ when $\xi$ decreases from $+\i$. Clearly, $\frac{d\Theta}{d\xi}=0$ and $\frac{dU}{d\xi}\leq0$ on $l_4$, any integral curve of $\eqref{keyequation}_1$ or \eqref{change} passing through the points on $l_4$ has the tangent vector pointing horizontally to the left. Likewise, any integral curve of $\eqref{keyequation}_1$ or \eqref{change} passing through the points on $l_5$ has the tangent vector pointing straightly upwards. Therefore, the curve $\Gamma_2$ can not intersect with $l_4$ or $l_5$. Moreover, $\Gamma_2$ can not terminate inside Region II when $\xi$ decreases from $+\i$ due to the fact that there is no equilibrium point to $\eqref{keyequation}_1$ or \eqref{change} inside Region II.  Therefore, the curve $\Gamma_2$ must tend to the equilibrium point $S_2\left(\a_1u_+,\a_2\t_+ \right)$ as $\xi\to-\i$. When the boundary data $(u_-,\t_-)\in\Gamma_2$ and $\frac{u_-}{v_-}=\frac{u_+}{v_+}$, the boundary layer solution to \eqref{BL-equation} exists if $\sqrt{\frac{\g-1}{2\g}}< M_+<1$.
	
	(b). $0<M_+\leq\sqrt{\frac{\g-1}{2\g}}$. In this case, the equilibrium point $S_2\left(\a_1u_+,\a_2\t_+ \right)$ is not above the horizontal line $\t=0$. Similar to the discussions in the previous case (a), the curve $\Gamma_2$ can be extended inside Region II until a unique point $Z_2$ on the line $\t=0$. Therefore, the curve $\Gamma_2$ lines between the point $Z_2$ on the line $\t=0$ and the equilibrium point $S_1(u_+, \t_+)$.
	
	Therefore,  the boundary layer solution to \eqref{BL-equation} exists if and only if $(u_-,\t_-)\in \Gamma_1\cup\Gamma_2$ and $\frac{u_-}{v_-}=\frac{u_+}{v_+}$ in Case III  (Subsonic case).
	
	The decay rates \eqref{trans-decay} and \eqref{subs-decay} follow the exactly same arguments as in \cite{Qin-Wang-2009}, and we omit the proof details for simplicity. Hence we complete the proof of Theorem \ref{maintheorem}.

	\section{Further discussions}
	\setcounter{equation}{0}

	\subsection{Stability of large boundary layer solution}
	
	Theorem \ref{maintheorem} gives the sufficient and necessary conditions on the existence/non-existence of large-amplitude boundary layer solutions \eqref{BL-equation} to the inflow problem  \eqref{(1.1)}--\eqref{initial}, $\eqref{boundarycases}_1$ or \eqref{Lag-equation}. 
	The next natural question is the time-asymptotic stability these large-amplitude boundary layer solutions. Huang, Li and Shi \cite{Huang-Li-Shi}, Qin and Wang \cite{Qin-Wang-2009, Qin-Wang-2011}, and Nakamura and Nishibata \cite{Nakamura-Nishibata}  proved the time-asymptotic stability of these boundary layer solutions and even composited with other basic wave patterns by the energy method and the Poincar\'e-type techniques for $\mathbb{R}_+$ in \cite{Nikkuni-Kawashima}, under the crucial smallness conditions on both the  wave amplitude and the initial perturbation, while similar stability results with large initial perturbation were later obtained by Bian, Fan, He and Zhao \cite{Bian-Fan-He-Zhao}, Fan, Liu, Wang and Zhao \cite{Fan-Liu-Wang-Zhao}, and Hong and Wang \cite{Hong-Wang-2017-1, Hong-Wang-2017-2}. Nevertheless, all the results mentioned above rely on the smallness of wave amplitude, which can not be extended to the large-amplitude case directly. In fact, it is still largely open to prove rigorously the time-asymptotic stability of these large-amplitude boundary layer solutions in Theorem \ref{maintheorem}. Some new insights and methods need to be developed.

	\subsection{Large boundary layer solutions to outflow problem}
	
	For the outflow problem  \eqref{(1.1)}--\eqref{initial}, $\eqref{boundarycases}_3$ with $u_-<0$, Kawashima, Nakamura, Nishibata and Zhu  \cite{KNNZ} established the existence of small-amplitude boundary layer solutions, when the right end state $(\r_+, u_+, \t_+)$ belongs to the supersonic, transonic and subsonic region respectively. However, compared with the inflow case explicitly characterized in Theorem \ref{maintheorem},  deriving the sufficient and necessary conditions for the existence/non-existence of large-amplitude boundary layer solutions to the outflow problem is more subtle and faces two new difficulties:
	
	\begin{itemize}
		
		\item ‌The integral curves of $\eqref{keyequation}_1$ may lose the monotonicity;
		
		\item The integral curves of $\eqref{keyequation}_1$ may be unbounded in $(u,\t)$-space.

	\end{itemize}
	
	Consequently, it is hard to trace the trajectories for the outflow problem and we leave the existence of its large boundary layer solutions for the future study.


\small
	\begin{thebibliography}{99}
		
		\bibitem{Bian-Fan-He-Zhao}
		\newblock   D. Bian, L. Fan, L. He, and H. Zhao,
		\newblock \emph{Viscous shock wave to an inflow problem for compressible viscous gas with large density oscillations},
		\newblock Acta Math. Appl. Sin. Engl. Ser., \textbf{35} (2019), 129--157.	
		
		\bibitem{Bressan-2007}
		\newblock   A. Bressan,
		\newblock \emph{Tutorial on the center manifold theorem},
		\newblock Lecture Notes in Math., 2007.
		
		\bibitem{Carr}
		\newblock J. Carr,
		\newblock ``Applications of centre manifold theory",
		\newblock Applied Mathematical Sciences, 35, Springer, New York-Berlin, 1981.
		
		\bibitem{Fan-Liu-Wang-Zhao}
		\newblock L. Fan, H. Liu, T. Wang, and H. Zhao,
		\newblock \emph{Inflow problem for the one-dimensional compressible Navier-Stokes equations under large initial perturbation},
		\newblock J. Differential Equations, \textbf{257} (2014), 3521--3553.	
		
		
		
		
		
		
		\bibitem{FD}
		\newblock H. Freistühler and D. Serre,
		\newblock \emph{$L^1$-stability of shock waves in scalar viscous conservation laws},
		\newblock Comm. Pure Appl. Math., \textbf{51} (1998), 291--301.	
		
		
		\bibitem{Gilbarg}
		\newblock D. Gilbarg,
		\newblock \emph{The existence and limit behavior of the one-dimensional shock layer},
		\newblock Amer. J. Math., \textbf{73} (1951), 256--274.
		
		
		\bibitem{Goodman}
		\newblock J. Goodman,
		\newblock \emph{Nonlinear asymptotic stability of viscous shock profiles for conservation laws},
		\newblock Arch. Ration. Mech. Anal., \textbf{95} (1986), 325--344.
		
		\bibitem{HKK}
		\newblock S. Han, M.-J. Kang, and J. Kim,
		\newblock \emph{Large-time behavior of composite waves of viscous shocks for the barotropic Navier-Stokes equations},
		\newblock SIAM J. Math. Anal., \textbf{55} (2023), 5526--5574.
		
		
		\bibitem{Hong-Wang-2017-1}
		\newblock H. Hong and T. Wang,
		\newblock \emph{Stability of stationary solutions to the inflow problem for full compressible Navier-Stokes equations with a large initial perturbation},
		\newblock SIAM J. Math. Anal., \textbf{49} (2017), 2138--2166.
		
		
		\bibitem{Hong-Wang-2017-2}
		\newblock H. Hong and T. Wang,
		\newblock \emph{Large-time behavior of solutions to the inflow problem of full compressible Navier-Stokes equations with large perturbation},
		\newblock Nonlinearity, \textbf{30} (2017), 3010--3039.
		
		
		
		\bibitem{Huang-Li-Matsumura}
		\newblock F. Huang, J. Li, and A. Matsumura,
		\newblock \emph{Asymptotic stability of combination of viscous contact wave with rarefaction waves for one-dimensional compressible Navier-Stokes system},
		\newblock Arch. Ration. Mech. Anal., \textbf{197} (2010), 89--116.
		
		
		
		\bibitem{Huang-Li-Shi}
		\newblock F. Huang, J. Li, and X. Shi,
		\newblock \emph{Asymptotic behavior of solutions to the full compressible Navier-Stokes equations in the half space},
		\newblock Commun. Math. Sci., \textbf{8} (2010), 639--654.
		
		
		\bibitem{Huang-Matsumura-Shi}
		\newblock F. Huang, A. Matsumura, and X. Shi,
		\newblock \emph{Viscous shock wave and boundary layer solution to an inflow problem for compressible viscous gas},
		\newblock Comm. Math. Phys., \textbf{239} (2003), 261--285.
		
		
		\bibitem{Huang-Qin}
		\newblock F. Huang and X. Qin,
		\newblock \emph{Stability of boundary layer and rarefaction wave to an outflow problem for compressible Navier-Stokes equations under large perturbation},
		\newblock J. Differential Equations, \textbf{246} (2009), 4077--4096.
		
		
		\bibitem{Huang-Xin-Yang}
		\newblock F. Huang, Z. Xin, and T. Yang,
		\newblock \emph{Contact discontinuities with general perturbation for gas motion},
		\newblock Adv. Math., \textbf{219} (2008), 1246--1297.
		
		
		\bibitem{Il'in-Oleinik} 
		\newblock A. M. Il'in and O. A. Oleinik,
		\newblock \emph{Asymptotic behavior of solution of the Cauchy problem for some quasi-linear equations for large values of the time},
		\newblock Mat. Sb. (N.S.), \textbf{51(93)} (1960), 191--216.
		
		
		
		\bibitem{KV}
		\newblock M.-J. Kang and A. Vasseur,
		\newblock \emph{{$L^2$}-contraction for shock waves of scalar viscous conservation laws},
		\newblock Ann. Inst. H. Poincar\'e{} C Anal. Non Lin\'eaire,  \textbf{34} (2017), 139--156.
		
		\bibitem{Kang-Vasseur-Wang-2023}
		\newblock M.-J. Kang, A. Vasseur, and Y. Wang,
		\newblock \emph{Time-asymptotic stability of composite waves of viscous shock and rarefaction for barotropic Navier-Stokes equations},
		\newblock Adv. Math., \textbf{419} (2023), Paper No. 108963, 66 pages.
		
		
		
		\bibitem{Kang-Vasseur-Wang-2024}
		\newblock M.-J. Kang, A. Vasseur, and Y. Wang,
		\newblock \emph{Time-asymptotic stability of generic Riemann solutions for compressible Navier-Stokes-Fourier equations},
		\newblock arxiv, (2023).
		
		\bibitem{KNNZ}
		\newblock S. Kawashima, T. Nakamura, S. Nishibata, and P. Zhu,
		\newblock \emph{Stationary waves to viscous heat-conductive gases in half-space: existence, stability and convergence rate},
		\newblock Math. Models Methods Appl. Sci., \textbf{20} (2010), 2201--2235.
		
		
		
		
		
		\bibitem{KNZ}
		\newblock S. Kawashima, S. Nishibata, and P. Zhu,
		\newblock \emph{Asymptotic stability of the stationary solution to the compressible Navier-Stokes equations in the half space},
		\newblock Comm. Math. Phys., \textbf{240} (2003), 483--500.
		
		
		\bibitem{KZ-08}
		\newblock S. Kawashima and P. Zhu,
		\newblock \emph{Asymptotic stability of nonlinear wave for the compressible Navier-Stokes equations in the half space},
		\newblock J. Differential Equations, \textbf{244} (2008), 3151--3179.
		
		\bibitem{KZ-09}
		\newblock S. Kawashima and P. Zhu,
		\newblock \emph{Asymptotic stability of rarefaction wave for the Navier-Stokes
			equations for a compressible fluid in the half space},
		\newblock Arch. Ration. Mech. Anal., \textbf{194} (2009), 105--132.
		
		
		
		\bibitem{Liu}
		\newblock T. Liu,
		\newblock \emph{Nonlinear stability of shock waves for viscous conservation laws},
		\newblock Mem. Amer. Math. Soc.,  \textbf{56} (1985), v+108.
		
		
		
		\bibitem{LX}
		\newblock T. Liu and Z. Xin,
		\newblock \emph{Nonlinear stability of rarefaction waves for compressible Navier-Stokes equations}
		\newblock Comm. Math. Phys.,  \textbf{118} (1988), 451--465.
		
		
		
		\bibitem{LZ}
		\newblock T. Liu and Y. Zeng,
		\newblock \emph{Shock waves in conservation laws with physical viscosity},
		\newblock Mem. Amer. Math. Soc.,  \textbf{234} (2015), vi+168.
		
		
		\bibitem{MZ}
		\newblock C. Mascia and K. Zumbrun,
		\newblock \emph{Stability of small-amplitude shock profiles of symmetric hyperbolic-parabolic systems},
		\newblock Comm. Pure Appl. Math.,  \textbf{57} (2004), 841--876.
		
		\bibitem{Matsumura-2001}
		\newblock A. Matsumura,
		\newblock \emph{Inflow and outflow problems in the half space for a one-dimensional isentropic model system of compressible viscous gas},
		\newblock  Methods Appl. Anal., \textbf{8} (2001), 645--666.
		
		
		
		\bibitem{MN-1985}
		\newblock A. Matsumura and K. Nishihara,
		\newblock \emph{On the stability of traveling wave solutions of a one-dimensional model system for compressible viscous gas},
		\newblock  Japan J. Appl. Math., \textbf{2} (1985), 17--25.
		
		
		\bibitem{MN-1986}
		\newblock A. Matsumura and K. Nishihara,
		\newblock \emph{Asymptotics toward the rarefaction wave of the solutions of a one-dimensional model system for compressible viscous gas},
		\newblock  Japan J. Appl. Math., \textbf{3} (1986), 1--13.
		
		
		
		
		
		
		
		\bibitem{Matsumura-Nishihara-2001}
		\newblock A. Matsumura and K. Nishihara,
		\newblock \emph{Large-time behaviors of solutions to an inflow problem in the half space for a one-dimensional system of compressible viscous gas},
		\newblock  Comm. Math. Phys., \textbf{222} (2001), 449--474.
		
		
		
		\bibitem{Nakamura-Nishibata}
		\newblock T. Nakamura and S. Nishibata,
		\newblock \emph{Stationary wave associated with an inflow problem in the half line for viscous heat-conductive gas},
		\newblock  J. Hyperbolic Differ. Equ., \textbf{8} (2011), 651--670.
		
		
		
		
		
		
		
		
		
		\bibitem{Nikkuni-Kawashima}
		\newblock Y. Nikkuni and S. Kawashima,
		\newblock \emph{Stability of stationary solutions to the half-space problem  for the discrete Boltzmann equation with multiple collisions},
		\newblock Kyushu J. Math., \textbf{54} (2000), 233--255.
		
		
		\bibitem{Nishihara-Yang-Zhao}
		\newblock K. Nishihara, T. Yang, and H. Zhao,
		\newblock \emph{Nonlinear stability of strong rarefaction waves for compressible
			Navier-Stokes equations},
		\newblock SIAM J. Math. Anal.,  \textbf{35} (2004), 1561--1597.
		
		
		
		\bibitem{Qin-Wang-2009}
		\newblock X. Qin and Y. Wang,
		\newblock \emph{Stability of wave patterns to the inflow problem of full compressible Navier-Stokes equations},
		\newblock SIAM J. Math. Anal., \textbf{41} (2009), 2057--2087.
		
		\bibitem{Qin-Wang-2011}
		\newblock X. Qin and Y. Wang,
		\newblock \emph{Large-time behavior of solutions to the inflow problem of full compressible Navier-Stokes equations},
		\newblock SIAM J. Math. Anal., \textbf{43} (2011), 341--366.
		
		
		\bibitem{Smoller} 
		\newblock J. Smoller,
		\newblock ``Shock Waves and Reaction-Diffusion Equations," New York: Springer, 1994.
		
		
		\bibitem{Szepessy-Xin}
		\newblock A. Szepessy and Z. Xin,
		\newblock \emph{Nonlinear stability of viscous shock waves},
		\newblock Arch. Ration. Mech. Anal.,  \textbf{122} (1993), 53--103.
		
		
		\bibitem{Zhang-Ding-Huang-Dong}
		\newblock   Z. Zhang, T. Ding, W. Huang, and Z. Dong,
		\newblock ``Qualitative Theory of Differential Equations" (in Chinese),
		\newblock Science Press, 1985.
		
	\end{thebibliography}
\end{document}